\documentclass[11pt,oneside,english]{amsart}
\usepackage[T1]{fontenc}
\usepackage{geometry}
\usepackage{amssymb}

\makeatletter
 \theoremstyle{plain}    
 \newtheorem{thm}{Theorem}[section]
 \numberwithin{equation}{section} 
 \numberwithin{figure}{section} 
 \theoremstyle{plain}
 \theoremstyle{remark}
 \newtheorem{rem}[thm]{Remark}
 \theoremstyle{plain}    
 \newtheorem{prop}[thm]{Proposition} 
 \theoremstyle{definition}
  \newtheorem{example}[thm]{Example}

\newcommand\rest{\hbox{\raise.17em\hbox{$ |\kern-.2em$}\lower.23em\hbox{$-$}}}
\newcommand\avint{\hbox{\hbox{$\displaystyle \int$}\hbox{\kern-.9em{$-$}}}}
\newcommand\smavint{\hbox{\hbox{$\int$}\hbox{\kern-.75em{$-$}}}}
\newcommand{\R}{\mathbb{R}}
\newcommand{\C}{\mathbb{C}}
\newcommand{\Z}{\mathbb{Z}}

\usepackage{babel}
\makeatother
\begin{document}

\title{Chern-Simons forms on associated bundles, and boundary terms}

\author{David L. Johnson}

\date{\today}

\address{Department of Mathematics\\
Lehigh University\\
Bethlehem, Pennsylvania 18015-3174}

\email{david.johnson@lehigh.edu}



\begin{abstract}
Let $E$ be a principle bundle over a compact manifold $M$ with compact
structural group $G$. For any $G$-invariant polynomial $P$, The
transgressive forms $TP(\omega)$ defined by Chern and Simons in \cite{CS}
are shown to extend to forms $\Phi P(\omega)$ on associated bundles
$B$ with fiber a quotient $F=G/H$ of the group. These forms satisfy
a heterotic formula \[
d\Phi P(\omega)=P(\Omega)-P(\Psi),\]
relating the characteristic form $P(\Omega)$ to a fiber-curvature
characteristic form. For certain natural bundles $B$, $P(\Psi)=0$,
giving a true transgressive form on the associated bundle, which leads
to the standard obstruction properties of characteristic classes as
well as natural expressions for boundary terms.
\end{abstract}
\maketitle

\section*{Introduction }

In their groundbreaking paper \cite{CS}, S-S Chern and James Simons
explain that their theory of what are now known as Chern-Simons classes

\begin{quotation}
grew out of an attempt to derive a purely combinatorial formula for
the first Pontryagin number of a 4-manifold. ... This process got
stuck by the emergence of a boundary term which did not yield to a
simple combinatorial analysis. The boundary term seemed interesting
in its own right and ....
\end{quotation}
Their {}``boundary term'' was in fact a geometric realization of
the transgression co-chains which occur in the Leray-Serre spectral
sequence of a principal bundle \cite{Borel}, and their importance
has grown out of the fact that, on the base manifold, they measure
finer geometric information than the primary characteristic classes
of the bundle. Parts of these transgressive forms have recently re-appeared
very much as boundary terms in a more classical sense \cite{Morgan-Johnson}.
Interestingly, one of the results established below re-constructs
forms in the unit tangent bundle which were originally constructed
by Chern in 1943 in conjunction with his version of the generalized
Gauss-Bonnet theorem \cite{Chern}. That construction of Chern is
interpreted in terms of his and Simons' later work in a broader context,
which also realizes components of these secondary characteristic classes
within more general bundles.

\section{Intermediate bundles and Chern-Simons forms}

Let $M$ be a compact, $n$-dimensional manifold, and let $\pi:E\rightarrow M$
be a principal bundle over $M$ with compact structural group $G$.
Real characteristic classes of $E$ can be determined by forms $P(\Omega)$,
where $P\in\mathcal{I}(\frak{g})$ is an adjoint-invariant polynomial
on the Lie algebra $\frak{g}$ of $G$, and $\Omega$ is the curvature
form of a connection $\omega$ on $E$. Such forms $P(\Omega)$ are
horizontal, invariant forms on $E$, so are naturally defined as forms
on $M$ itself. 

Invariance properties of these polynomial forms, along with the Bianchi
identity, traditionally are used to show that these forms are closed
on $M$ (cf. \cite{Milnor-Stasheff}). Moreover, the forms $P(\Omega)$
in the cohomology of $E$ itself are exact, which not only verifies
that $P(\Omega)$ are closed on $M$ but also implies the existence
of forms $TP(\omega)$ on $E$, primitives of the characteristic forms
$P(\Omega)$. Given a connection $\omega$ on $E$, Chern and Simons
derive in \cite{CS} an explicit formula for these transgressive forms.

\begin{thm}
{[}Chern-Simons{]}. Let $\pi:E\rightarrow M$ be a principal bundle
over a compact $n$-manifold $M$ with compact structural group $G$.
If $P\in\mathcal{I}^{k}(\frak{\frak{g}})$ is a degree-$k$, adjoint-invariant
polynomial on the Lie algebra $\frak{g}$ of $G$, and if $\omega$
is a connection on $E$, then the $(2k-1)$-form \[
TP(\omega):=\sum_{i=0}^{k-1}A_{i}P(\omega,[\omega,\omega]^{i},\Omega^{k-i-1})\]
is a $G$-invariant form on $E$ satisfying $dTP(\omega)=P(\Omega)$,
where $A_{i}:=(-1)^{i}k!(k-1)!/2^{i}(k+1)!(k-1-i)!$ and $P$ is realized
as a symmetric, multilinear functional $P:\frak{g}\times\cdots\times\frak{g}\rightarrow\R$
by polarization. 
\end{thm}
Let now $\pi_{2}:B\rightarrow M$ be an associated bundle to the principal
bundle $\pi:E\rightarrow M$ as before, with fibers $F$ which are
homogeneous spaces, quotients of the structural group $G$ by the
isotropy subgroup $H$ of the right action of $G$ on $F$. Such a
bundle will be called an \emph{intermediate bundle} of $E$. The primary
example of this situation is when $E$ is the bundle of oriented frames
of an oriented Riemannian manifold $M$, and $B$ the unit tangent
bundle of $M$. In that case $G=SO(n)$ and $H=SO(n-1)$. 

The two bundles are related, and in fact the total space of $E$ is
a principal bundle with group $H$ over $B$: \begin{equation}
\begin{array}{ccc}
E & \begin{array}{c}
\pi_{1}\\
\longrightarrow\\
\,\end{array} & B\\
\begin{array}{ccc}
\pi & \big{\downarrow} & \,\,\end{array} &  & \begin{array}{ccc}
\pi_{2} & \big{\downarrow} & \,\,\,\end{array}\\
M & \begin{array}{c}
\,\\
=\\
\,\end{array} & M\end{array}\label{eqn:intermediate-bundle-diagram}\end{equation}
Decompose the connection $\omega$ in terms of the $G$-equivariant
distribution $\overline{\frak{h}}:=ker((\pi_{1})_{*})$ on $E$, $\omega=\phi+\psi$,
where $\psi=\left.\omega\right|_{ker((\pi_{1})_{*})}$. $\psi$ is
the connection induced from $\omega$ on the principal bundle $\pi_{1}:E\rightarrow B$
with structural group $H$. Since $H$ is a reductive subgroup, if
(by choice of bases) at $p\in E$ $\left.\omega\right|_{ker((\pi_{1})_{*})}:ker((\pi_{1})_{*})\rightarrow\frak{h}$,
then $\phi$ takes values in a reductive complement $\frak{p}$ to
$\frak{h}$. Also, $[\psi,\psi]$ has values in $\frak{h}$, and $[\psi,\phi]$
takes values in $\frak{p}$.

If $\mathcal{H}$ represents the $\psi$-horizontal projection (and
subspace), the curvature of the two connections are related by $\Psi=d\psi+\frac{1}{2}[\psi,\psi]$
and \begin{eqnarray*}
\Omega & = & d\omega+\frac{1}{2}[\omega,\omega]\\
 & = & d\phi+[\psi,\phi]+\frac{1}{2}[\phi,\phi]+\Psi\\
 & = & d_{\mathcal{H}}\phi+\Psi+\frac{1}{2}[\phi,\phi].\end{eqnarray*}
The $\psi$-covariant differential $d_{\mathcal{H}}\phi:=d\phi+[\psi,\phi]$
is the restriction of $d\phi$ to $\psi$-horizontal tangents. Similarly,
by the Bianchi identity, \begin{eqnarray*}
d_{\mathcal{H}}\Omega & = & d\Omega+[\psi,\Omega]\\
 & = & [\Omega,\omega]+[\psi,\Omega]\\
 & = & [\Omega,\omega]-[\Omega,\psi]\\
 & = & [\Omega,\phi].\end{eqnarray*}

\begin{thm}
\label{thm:Psi-P-omega} Let $\pi:E\rightarrow M$ be a principal
bundle over a compact $n$-manifold $M$ with compact structural group
$G$. Let $P\in\mathcal{I}^{k}(\frak{g})$ be a degree-$k$, adjoint-invariant
polynomial on the Lie algebra $\frak{g}$ of $G$, and let $\omega$
be a connection on $E$. For an intermediate bundle $\pi_{2}:B\rightarrow M$
with fiber $G/H$ as above, the form \[
\Phi P(\omega):=\sum_{i=0}^{k-1}\sum_{j=0}^{k-1-i}A_{ij}P(\phi,[\phi,\phi]^{i},\Psi^{j},\Omega^{k-i-j-1}),\]
 where $A_{ij}:=(-1)^{i}\frac{(i+j)!(k-j-1)!k!}{2^{i}(k-i-j-1)!i!(k+i)!j!}$,
is a $\pi_{1}$-horizontal, $Ad_{H}$-invariant form on $E$, thus
represents a form on $B$. In addition, \begin{equation}
d\Phi P(\omega)=P(\Omega)-P(\Psi).\label{eqn:d-T-Phi}\end{equation}
 
\end{thm}
\begin{proof}
That $\Phi P(\omega)$ is $\pi_{1}$-horizontal and $Ad_{H}$-invariant,
and so is a form on $B$, follows by the definitions of $\phi$ and
$\Psi$, and the fact that $P$ is invariant under $Ad_{G}$. Also,
for any invariant polynomial $P$ and equivariant, $\psi$-horizontal
forms $\alpha_{1},\ldots,\alpha_{k}$ of degrees $p_{1},\ldots,p_{k}$,
respectively, it is straightforward (cf. \cite{CS}) that, on $B$,
\begin{eqnarray*}
dP(\alpha_{1},\ldots,\alpha_{k}) & = & \sum_{i}(-1)^{p_{1}+\cdots+p_{i-1}}P(\alpha_{1},\ldots,\alpha_{i-1},d_{\mathcal{H}}\alpha_{i},\alpha_{i+1},\ldots,\alpha_{k}).\end{eqnarray*}
 We now show that there are constants $A_{ij}$ satisfying \ref{eqn:d-T-Phi},
and that they are well-defined. For any choices of $A_{ij}$, \begin{eqnarray*}
d\Phi P(\omega) & = & \sum_{i=0}^{k-1}\sum_{j=0}^{k-1-i}A_{ij}dP(\phi,[\phi,\phi]^{i},\Psi^{j},\Omega^{k-i-j-1})\\
 & = & \sum_{i=0}^{k-1}\sum_{j=0}^{k-1-i}A_{ij}P(d_{\mathcal{H}}\phi,[\phi,\phi]^{i},\Psi^{j},\Omega^{k-i-j-1})\\
 &  & -iA_{ij}P(\phi,2[d_{\mathcal{H}}\phi,\phi],[\phi,\phi]^{i-1},\Psi^{j},\Omega^{k-i-j-1})\\
 &  & -jA_{ij}P(\phi,[\phi,\phi]^{i},d_{\mathcal{H}}\Psi,\Psi^{j-1},\Omega^{k-i-j-1})\\
 &  & -(k-i-j-1)A_{ij}P(\phi,[\phi,\phi]^{i},\Psi^{j},d_{\mathcal{H}}\Omega,\Omega^{k-i-j-2})\\
 & = & \sum_{i=0}^{k-1}\sum_{j=0}^{k-1-i}A_{ij}P(\Omega-\Psi-\frac{1}{2}[\phi,\phi],[\phi,\phi]^{i},\Psi^{j},\Omega^{k-i-j-1})\\
 &  & -iA_{ij}P(\phi,2[\Omega-\Psi-\frac{1}{2}[\phi,\phi],\phi],[\phi,\phi]^{i-1},\Psi^{j},\Omega^{k-i-j-1})\\
 &  & -jA_{ij}P(\phi,[\phi,\phi]^{i},0,\Psi^{j-1},\Omega^{k-i-j-1})\\
 &  & -(k-i-j-1)A_{ij}P(\phi,[\phi,\phi]^{i},\Psi^{j},[\Omega,\phi],\Omega^{k-i-j-2})\\
 & = & \sum_{i=0}^{k-1}\sum_{j=0}^{k-1-i}A_{ij}P([\phi,\phi]^{i},\Psi^{j},\Omega^{k-i-j})\\
 &  & -A_{ij}P([\phi,\phi]^{i},\Psi^{j+1},\Omega^{k-i-j-1})-\frac{1}{2}A_{ij}P([\phi,\phi]^{i+1},\Psi^{j},\Omega^{k-i-j-1})\\
 &  & -2iA_{ij}P(\phi,[\Omega,\phi],[\phi,\phi]^{i-1},\Psi^{j},\Omega^{k-i-j-1})\\
 &  & +2iA_{ij}P(\phi,[\Psi,\phi],[\phi,\phi]^{i-1},\Psi^{j},\Omega^{k-i-j-1})\\
 &  & -(k-i-j-1)A_{ij}P(\phi,[\phi,\phi]^{i},\Psi^{j},[\Omega,\phi],\Omega^{k-i-j-2}).\end{eqnarray*}
Re-grouping by the powers of $[\phi,\phi]$ and $\Psi$, and using
the identity from \cite{CS} which comes from invariance of the polynomial,
\[
0=(-1)^{p_{1}}P([\alpha_{1},\phi],\alpha_{2},\ldots,\alpha_{k})+\cdots+(-1)^{p_{1}+\cdots+p_{k}}P(\alpha_{1},\ldots,[\alpha_{k},\phi]),\]
 so that, in particular, \begin{eqnarray*}
P([\phi,\phi]^{i},\Psi^{j},\Omega^{k-i-j}) & = & -jP(\phi,[\phi,\phi]^{i-1},[\Psi,\phi],\Psi^{j-1},\Omega^{k-i-j})\\
 &  & -(k-i-j)P(\phi,[\phi,\phi]^{i-1},\Psi^{j},[\Omega,\phi],\Omega^{k-i-j-1})\end{eqnarray*}
for $\phi$ a $\frak{g}$-valued 1-form, (setting $A_{i,j}=0$ if
either $i$ or $j$ is negative, or if $i+j>k-1$) \begin{eqnarray*}
d\Phi P(\omega) & = & \sum_{i=0}^{k-1}\sum_{j=0}^{k-1-i}(A_{ij}-A_{i,j-1}-\frac{1}{2}A_{i-1,j})P([\phi,\phi]^{i},\Psi^{j},\Omega^{k-i-j})\\
 &  & -(2iA_{ij}+(k-i-j)A_{i-1,j})P(\phi,[\phi,\phi]^{i-1},\Psi^{j},[\Omega,\phi],\Omega^{k-i-j-1})\\
 &  & +2iA_{i,j-1}P(\phi,[\phi,\phi]^{i-1},[\Psi,\phi],\Psi^{j-1},\Omega^{k-i-j})\\
 & = & \sum_{i=0}^{k-1}\sum_{j=0}^{k-1-i}(A_{ij}-A_{i,j-1}-\frac{1}{2}A_{i-1,j})P([\phi,\phi]^{i},\Psi^{j},\Omega^{k-i-j})\\
 &  & -(2iA_{ij}+(k-i-j)A_{i-1,j})P(\phi,[\phi,\phi]^{i-1},\Psi^{j},[\Omega,\phi],\Omega^{k-i-j-1})\\
 &  & -2iA_{i,j-1}\frac{1}{j}\left(P([\phi,\phi]^{i},\Psi^{j},\Omega^{k-i-j})+(k-i-j)P(\phi,[\phi,\phi]^{i-1},\Psi^{j},[\Omega,\phi],\Omega^{k-i-j-1})\right)\\
 & = & \sum_{i=0}^{k-1}\sum_{j=0}^{k-1-i}(A_{ij}-A_{i,j-1}-\frac{1}{2}A_{i-1,j}-2\frac{i}{j}A_{i,j-1})P([\phi,\phi]^{i},\Psi^{j},\Omega^{k-i-j})\\
 &  & -(2iA_{ij}+(k-i-j)A_{i-1,j}+2\frac{i(k-i-j)}{j}A_{i,j-1})P(\phi,[\phi,\phi]^{i-1},\Psi^{j},[\Omega,\phi],\Omega^{k-i-j-1}).\end{eqnarray*}
We will have $d\Phi P(\omega)=P(\Omega)-P(\Psi)$ if the coefficient
of $P([\phi,\phi]^{i},\Psi^{j},\Omega^{k-i-j})$ is 1 for $i=j=0$
and -1 for $i=0,\, j=k$, and $0$ otherwise, as well as the coefficients
of $P(\phi,[\phi,\phi]^{i-1},\Psi^{j},[\Omega,\phi],\Omega^{k-i-j-1})$
vanishing. That is, for $(i,j)\neq(0,0)$, \begin{eqnarray*}
0 & = & A_{ij}-(\frac{j+2i}{j})A_{i,j-1}-\frac{1}{2}A_{i-1,j}\\
0 & = & 2iA_{ij}+(k-i-j)A_{i-1,j}+2\frac{i(k-i-j)}{j}A_{i,j-1}\end{eqnarray*}
or\begin{eqnarray*}
A_{i,j-1} & = & -\frac{j(k-j)}{2i(k+i)}A_{i-1,j}\\
A_{ij} & = & \left(\frac{(i+j)(i+j-k)}{2i(k+i)}\right)A_{i-1,j}.\end{eqnarray*}
There is a necessary consistency condition, in that two recursive
formulas be consistent, that is \[
A_{i+1,j-1}=\left(\frac{(i+j)(i+j-k)}{2(i+1)(k+i+1)}\right)A_{i,j-1}=\left(\frac{(i+j)(i+j-k)}{2(i+1)(k+i+1)}\right)\left(-\frac{j(k-j)}{2i(k+i)}\right)A_{i-1,j}\]
versus\[
A_{i+1,j-1}=-\frac{j(k-j)}{2(i+1)(k+i+1)}A_{i,j}=-\frac{j(k-j)}{2(i+1)(k+i+1)}\left(\frac{(i+j)(i+j-k)}{2i(k+i)}\right)A_{i-1,j},\]
which indeed do give the same expression, so that the double recursion
defining $A_{ij}$ is consistent. 

From the second recursion, setting $A_{0,0}=1$, we obtain\begin{eqnarray*}
A_{i,0} & = & \left(\frac{(i-k)}{2(k+i)}\right)A_{i-1,0}\\
 & = & \left(\frac{(i-k)}{2(k+i)}\right)\left(\frac{(i-1-k)}{2(k+i-1)}\right)A_{i-2,0}\\
 & = & \frac{(-1)^{i}k!(k-1)!}{2^{i}(k-i-1)!(k+i)!},\end{eqnarray*}
exactly agreeing with the terms $A_{i}$ of \cite{CS}, as expected.
Now, using the first recursion,\begin{eqnarray*}
A_{i,j} & = & -\frac{2(i+1)(k+i+1)}{j(k-j)}A_{i+1,j-1}\\
 & = & (-1)^{2}\left(\frac{2(i+1)(k+i+1)}{j(k-j)}\right)\left(\frac{2(i+2)(k+i+2)}{(j-1)(k-j+1)}\right)A_{i+2,j-2}\\
 & = & (-1)^{2}\left(\frac{2^{2}(i+1)(i+2)(k+i+1)(k+i+2)}{j(j-1)(k-j)(k-j+1)}\right)A_{i+2,j-2}\\
 & = & (-1)^{j}\frac{2^{j}(i+j)!(k+i+j)!(k-j-1)!}{i!(k+i)!j!(k-1)!}A_{i+j,0}\\
 & = & (-1)^{j}\frac{2^{j}(i+j)!(k+i+j)!(k-j-1)!}{i!(k+i)!j!(k-1)!}\frac{(-1)^{i+j}k!(k-1)!}{2^{i+j}(k-i-j-1)!(k+i+j)!}\\
 & = & (-1)^{i}\frac{(i+j)!(k-j-1)!k!}{2^{i}(k-i-j-1)!i!(k+i)!j!},\end{eqnarray*}
which is of course the general term.

The nature of the recursion will guarantee that the coefficients of
$P([\phi,\phi]^{i},\Psi^{j},\Omega^{k-i-j})$ will be 0 except when
$i=j=0$, or $j=k$, and that the coefficient of $P(\Omega^{k})$
will be 1, because $A_{ij}=0$ if either $i$ or $j$ is negative.
Also, the coefficient of $P(\Psi^{k})$ will be $-A_{0,k-1}$ (recalling
that $A_{ij}=0$ if $i+j>k-1$). Now, \[
-A_{0,k-1}=-\frac{(k-1)!k!}{(k)!(k-1)!}=-1,\]
as claimed.
\end{proof}
The right-hand side of (\ref{eqn:d-T-Phi}) is not, unfortunately,
exactly the characteristic class $P(\Omega)$ that one might hope
for. Fortunately, though, in certain circumstances it can be shown
that $P(\Psi)=0$, for which bundles $B$ the form $\Phi P(\omega)$
will represent a secondary characteristic form of $P(\Omega)$ on
the associated bundle. This occurs in particular for the Gauss-Bonnet
integrand on the unit tangent bundle, which gives the connection between
the Chern-Simons class $Te(\omega)$ of the Riemannian connection
of an even-dimensional Riemannian manifold $M$ and the formulas for
the boundary term described by Chern in \cite{Chern}. 

It is not the case that $\Phi P(\omega)$ is the $\psi$-horizontal
part of $TP(\omega)$, which instead is only the terms in $\Phi P(\omega)$
with $j=0$. The additional terms, those involving the curvature $\Psi$
of $\psi$, can be expressed in terms of $\Omega$, $\phi$, as explained
below.

\begin{rem}
The formula (\ref{eqn:d-T-Phi}) of Theorem (\ref{thm:Psi-P-omega})
\[
d\Phi P(\omega)=P(\Omega)-P(\Psi)\]
is a general version of the heterotic formula $dH=Tr(F\wedge F)-Tr(R\wedge R)$
of \cite{Vafa-Witten}, in that, in the case of a tensorial bundle,
the curvature term $\Omega$ is related to the curvature of the base
manifold, and the curvature $\Psi$ is a curvature of the fibers.
\end{rem}

\subsection{$\Psi$}

The curvature form of the connection $\psi$ on $E\rightarrow B$
can easily be interpreted in terms of $\Omega$ and $\phi$, since
$d_{\mathcal{H}}\phi=\Omega-\Psi-\frac{1}{2}[\phi,\phi]$ and since
$\Psi$ has values in the subalgebra $\frak{h}\subset\frak{g}$, but
$\phi$, hence $d_{\mathcal{H}}\phi$, takes values in the reductive
complement $\frak{p}$ of $\frak{h}$. Thus, \[
\Psi=\Omega_{\frak{h}}+\frac{1}{2}[\phi,\phi]_{\frak{h}},\]
where $\Omega_{\frak{h}}$ is the projection onto the subalgebra $\frak{h}$.

\subsection{Bundles naturally associated to characteristic classes }

The bundle $B$ is \emph{naturally associated} with a particular invariant
polynomial $P$ (or the characteristic class $P(\Omega)$) if $P(\Psi)=0$.
In such a case, Theorem (\ref{thm:Psi-P-omega}) implies that $\Phi P(\omega)$
will be a transgressive form for $P(\Omega)$ in the same sense as
in the principal bundle case. Since $\Psi$ is the curvature of the
connection $\psi$ on the bundle $E\rightarrow B$, $B$ will be naturally
associated with any polynomial which vanishes completely on $\frak{h}\subset\frak{g}$.
Examples of such bundles are:

\begin{enumerate}
\item If $E$ is the bundle of oriented frames of an oriented, even rank,
real vector bundle $\xi^{2k}$ (with of course a chosen inner product),
then set $B$ to be the unit sphere bundle of $\xi$, and $H=1\times SO(2k-1)\subset SO(2k)$.
In that case, since the Euler polynomial $e$ (which is up to a constant
the Pfaffian, the square root of the determinant) vanishes on the
subalgebra $\frak{h}=0\oplus so(2k-1)$, then $e(\Psi)=0$, and $d\Phi e(\omega)=P(\Omega)$. 
\item For $E$ the bundle of unitary frames of a rank-$k$ hermitian vector
bundle $\xi$, a $U(k)$-principal bundle, then for $\frak{h}=0\oplus u(k-1)\subset u(k)$,
the determinant is 0, so $c_{k}(\Psi)=0$ if $B$ is the bundle with
fiber $U(k)/I\times U(k-1)$, which is of course again just the unit
$S^{2k-1}$-sphere bundle of $\xi$. 
\item For the same unitary frame bundle $E$, then taking $\frak{h}=SU(k)$,
so that $B$ is the $U(1)$-principal bundle of the determinant bundle
of $\xi$, then $c_{1}(\Psi)=0$ whenever $\Psi\in\frak{h}$, so (as
expected) the determinant bundle is naturally associated to $c_{1}$. 
\item Again, on the unitary frame bundle $E=F(\xi)$, taking $H=I_{k-j+1}\times U(j-1)$,
so that $B$ is the Stiefel bundle of $(k-j+1)$-frames on $\xi$,
then $c_{j}(\Psi)=0$ whenever $\Psi\in\frak{h}$, and this Stiefel
bundle is naturally associated to the Chern class $c_{j}(\Omega)$. 
\end{enumerate}

\section{Obstructions}

For specific intermediate bundles, the characteristic classes $P(\Omega)$
are obstructions to the existence of global sections. Using the forms
$\Phi P(\omega)$, following Chern's original construction, the characteristic
classes $P(\Omega)$ can be explicitly computed as obstructions. The
same proof, applied to integration over chains rather than cycles,
gives relative versions of each of these classes. The boundary term
will depend upon a choice of section on the boundary, so may not directly
lead to combinatorial formulas for the characteristic classes. In
the case of the Euler class of the tangent bundle, the normal field
of the boundary gives a canonical section of the tangent bundle over
the boundary.

\subsection{The Euler class}

This first result, in the case of a cycle, is of course classical,
and the method of proof is essentially that of \cite{Chern}. In the
general case, the result does follow from Chern's result, but was
not stated as such by him. Several authors have presented proofs of
the general result, usually just for the tangent bundle of a manifold-with-boundary,
such as \cite{Hutchings,Li-Duan,Spivak-V}. The formulations differ
from case to case, but each basically recovers Chern's transgression
form, as does the current version.

\begin{prop}
Let $\xi$ be a rank-$2k$, oriented vector bundle over a compact
manifold $M$, with a smooth metric. Let $\omega$ be a metric-compatible
connection on $\xi$. Let $\alpha$ be a smooth singular $2k$-chain
in $M$, and let $\sigma$ be a generic section of $\xi$. Then\begin{eqnarray*}
\int_{\alpha}e(\Omega) & = & \sum_{j=1}^{m}a_{j}+\int_{\partial\alpha}s^{*}(\Phi e(\omega)),\end{eqnarray*}
where $\{ p_{1},\ldots,p_{m}\}$ are the zeros of $\sigma$, with
$a_{j}$ the index of the zero of $\sigma$ at $p_{j}$, and $s:=\sigma/\left\Vert \sigma\right\Vert $.
\end{prop}
\begin{proof}
For an oriented, rank-$2k$ vector bundle $\xi$ over $M$, and for
$\omega$ a metric-compatible connection on $\xi$, the naturally-associated
intermediate bundle $B$ is of course the unit $(2k-1)$-sphere bundle
$S(\xi)$ of $\xi$. Within $S(\xi)$, $d\Phi e(\omega)=e(\Omega)$,
since $e(\Psi)=0$, $\Psi$ lying in $so(2k-1)$. 

Since $\sigma$ is generic, we can assume that the zero-section of
$\sigma$ will have intersection with $\alpha$ a finite collection
of points $\{ p_{1},\ldots,p_{m}\}$ in the interiors of $2k$-simplices
of $\alpha$, with nonzero, finite-degree singularities. Then, for
any $\epsilon>0$ sufficiently small, $s:=\sigma/\left\Vert \sigma\right\Vert $
defines a section over $Supp(\alpha)\backslash\{ B(p_{1},\epsilon)\cup\cdots\cup B(p_{m},\epsilon)\}$,
where $B(p_{1},\epsilon)$ is the $\epsilon$-ball within the appropriate
$2k$-simplex, and \begin{eqnarray*}
\int_{\alpha}e(\Omega) & = & \lim_{\epsilon\downarrow0}\int_{\alpha\backslash B(p_{1},\epsilon)\cup\cdots\cup B(p_{m},\epsilon)}s^{*}(d\Phi e(\omega))\\
 & = & \lim_{\epsilon\downarrow0}\sum_{j=1}^{m}\int_{s_{*}(\partial B(p_{j},\epsilon))}\Phi e(\omega)+\int_{\partial\alpha}s^{*}(\Phi e(\omega)).\end{eqnarray*}
Since each singularity of $\sigma$ is of finite, nonzero degree $a_{j}$,
$\lim_{\epsilon\downarrow0}s_{*}(\partial B(p_{j},\epsilon))\sim a_{j}\pi^{-1}(p_{j})$
(homologous) where $\pi:S(\xi)\rightarrow M$ is the bundle projection.
Since $\Omega$ is $\pi:S(\xi)\rightarrow M$ horizontal, using the
form of $\Psi$ above and the fact that $[\phi,\phi]$ has image contained
in $\frak{h}$ in this case, \begin{eqnarray*}
\lim_{\epsilon\downarrow0}\sum_{j=1}^{m}\int_{s_{*}(\partial B(p_{j},\epsilon))}\Phi e(\omega) & = & \sum_{j=1}^{m}a_{j}\int_{\pi^{-1}(p_{j})}\Phi e(\omega)\\
 & = & \sum_{j=1}^{m}a_{j}\int_{\pi^{-1}(p_{j})}\sum_{i=0}^{k-1}A_{i,k-1-i}e(\phi,[\phi,\phi]^{i},\Psi^{k-1-i})\\
 & = & \sum_{j=1}^{m}a_{j}\int_{\pi^{-1}(p_{j})}\sum_{i=0}^{k-1}A_{i,k-1-i}\frac{1}{2^{k-1-i}}e(\phi,[\phi,\phi]^{k-1})\\
 & = & \sum_{j=1}^{m}a_{j}\int_{\pi^{-1}(p_{j})}\sum_{i=0}^{k-1}(-1)^{i}\frac{(k-1)!k!}{2^{k-1}(k+i)!(k-1-i)!}e(\phi,[\phi,\phi]^{k-1})\\
 & = & \sum_{j=1}^{m}a_{j}\int_{\pi^{-1}(p_{j})}\frac{k}{(2k-1)2^{k-1}}e(\phi,[\phi,\phi]^{k-1})\\
 & = & \sum_{j=1}^{m}a_{j},\end{eqnarray*}
since the integral can be normalized by applying it to the singularities
of the longitudinal flow on the sphere $S^{2k}$.
\end{proof}
\begin{rem}
If $\xi=T_{*}(M)$ and $M$ is a $2k$-manifold with boundary $\partial M$,
then the usual Gauss-Bonnet-Chern theorem, with boundary, can be recovered
by taking $\sigma$ to be the unit normal field to $\partial M\subset M$,
and of course the Poincar\'{e}-Hopf theorem. 

Using this result it will be possible to extend to all even dimensions
the usual method of proof of the 2-dimensional Gauss-Bonnet theorem,
integrating simplex by simplex. For each $2k$-simplex, take $\sigma$
to have a single, index-1 zero in the interior, and to be the unit
normal field on $\partial\sigma$. On $(2k-1)$-simplices in the boundary,
the boundary integrals will cancel due to the choice of normal. Careful
computation of the contributions from the lower-dimensional simplices
(taking limits of smooth approximations to the corners of $\sigma$)
will produce the result.
\end{rem}

\subsection{Chern Classes}

Since Chern classes are defined, by the splitting principle, from
the Euler class \cite{Milnor-Stasheff}, the situation is quite similar
for Chern classes as for the Euler class. For a complex rank-$k$
vector bundle $\xi\rightarrow M$, the transgression of the the $j^{th}$
Chern class $c_{j}(\xi)$ will be naturally-defined on the Stiefel
bundle $B:=V_{k-j+1}(\xi)$ of $(k-j+1)$-frames on $\xi$, with fiber
$U(k)/U(j-1)$. Within $V_{k-j+1}(\xi)$, $d\Phi c_{j}(\omega)=c_{j}(\Omega)$,
since $c_{j}(\Psi)=0$ for $\Psi$ lying in $u(j-1)$. 

\begin{prop}
Let $\xi$ be a rank-$k$, complex vector bundle over a compact manifold
$M$, with a smooth hermitian metric. Let $\omega$ be a metric-compatible
connection on $\xi$. Let $\alpha$ be a smooth singular $2j$-chain
in $M$, and let $(\sigma_{1},\ldots,\sigma_{k-j})$ be a unitary
$(k-j)$-frame of $\left.\xi\right|_{Supp(\alpha)}$. Let $\sigma$
be a generic section of $\xi/Span_{\mathbb{C}}\{\sigma_{1},\ldots,\sigma_{k-j}\}\cong\xi^{\perp}$
with no zeros on $\partial\alpha$. Then\begin{eqnarray*}
\int_{\alpha}c_{j}(\Omega) & = & \sum_{l=1}^{m}a_{l}+\int_{\partial\alpha}s^{*}\Phi c_{j}(\omega),\end{eqnarray*}
where $\{ p_{1},\ldots,p_{m}\}$ are the singularities of $s:=(\sigma_{1},\ldots,\sigma_{k-j},\sigma/\left\Vert \sigma\right\Vert )$
as a section of the Stiefel bundle $\left.V_{k-j+1}(\xi)\right|_{\alpha}$,
with $a_{l}$ the index of the singularity of $\sigma$ at $p_{l}$.
\end{prop}
\begin{proof}
Let $\alpha$ be a smooth singular $2j$-cycle in $M$. Since $rank_{\R}(\xi)>2j$,
there is a unitary $(k-j)$-frame $(\sigma_{1},\ldots,\sigma_{k-j})$
of $\left.\xi\right|_{Supp(\alpha)}$. Let $\sigma$ be a generic
section of the orthogonal complement $\xi^{\perp}\cong\xi/Span\{\sigma_{1},\ldots,\sigma_{k-j}\}$.
Since $\sigma$ is generic, we can assume that the zero-section of
$\sigma$ will have intersection with $Supp(\alpha)$ a finite collection
of points $\{ p_{1},\ldots,p_{m}\}$ in the interiors of $2j$-simplices
of $\alpha$, with nonzero, finite-degree singularities. Then, for
any $\epsilon>0$ sufficiently small, $s:=(\sigma_{1},\ldots,\sigma_{k-j+1})$,
with $\sigma_{k-j+1}:=\sigma/\left\Vert \sigma\right\Vert $, defines
a section of $V_{k-j+1}(\xi)$ over $Supp(\alpha)\backslash\{ B(p_{1},\epsilon)\cup\cdots\cup B(p_{m},\epsilon)\}$,
where $B(p_{1},\epsilon)$ is the $\epsilon$-ball within the appropriate
$2j$-simplex, and, since $d\Phi c_{j}(\omega)=c_{j}(\Omega)$ on
$B=V_{k-j+1}(\xi)$, \begin{eqnarray*}
\int_{\alpha}c_{j}(\Omega) & = & \lim_{\epsilon\downarrow0}\int_{\alpha\backslash B(p_{1},\epsilon)\cup\cdots\cup B(p_{m},\epsilon)}s^{*}(d\Phi c_{j}(\omega))\\
 & = & \lim_{\epsilon\downarrow0}\sum_{l=1}^{m}\int_{s_{*}(\partial B(p_{l},\epsilon))}\Phi c_{j}(\omega)+\int_{\partial\alpha}s^{*}\Phi c_{j}(\omega).\end{eqnarray*}
Since each singularity of $\sigma_{k-j+1}$ is of finite, nonzero
degree $a_{l}$ as a section of the unit sphere bundle in $\xi/Span\{\sigma_{1},\ldots,\sigma_{k-j}\}$,
$\lim_{\epsilon\downarrow0}s_{*}(\partial B(p_{l},\epsilon))=a_{l}S(p_{l})$
where $S(p_{l})$ is the $(2j-1)$-sphere in $\left.V_{k-j+1}(\xi)\right|_{p_{l}}$
defined by fixing $\sigma_{1},\ldots,\sigma_{k-j}$ at $p_{l}$, and
varying $\sigma_{k-j+1}(p_{l})$ among all unit vectors orthogonal
to the span of $\{\sigma_{1},\ldots,\sigma_{k-j}\}$. Since $\Omega$
is $\pi:V_{k-j+1}(\xi)\rightarrow M$ horizontal, using the form of
$\Psi$above, \begin{eqnarray*}
\lim_{\epsilon\downarrow0}\sum_{l=1}^{m}\int_{s_{*}(\partial B(p_{l},\epsilon))}\Phi c_{j}(\omega) & = & \sum_{l=1}^{m}a_{l}\int_{\pi^{-1}(p_{l})}\Phi c_{j}(\omega)\\
 & = & \sum_{l=1}^{m}a_{l}\int_{\pi^{-1}(p_{l})}\sum_{i=0}^{j-1}A_{i,j-1-i}c_{j}(\phi,[\phi,\phi]^{i},\Psi^{j-1-i})\\
 & = & \sum_{l=1}^{m}a_{l}\int_{\pi^{-1}(p_{l})}\sum_{i=0}^{k-1}A_{i,j-1-i}\frac{1}{2^{j-1-i}}c_{j}(\phi,[\phi,\phi]^{i},[\phi,\phi]_{\frak{h}}^{j-i-1})\\
 & = & \sum_{l=1}^{m}a_{l},\end{eqnarray*}
again by normalizing the integral on a test case, such as the sum
$\mathcal{O}(1)\oplus\cdots\oplus\mathcal{O}(1)$ of $j$ copies of
the line bundle with $c_{1}=1$ on $\C\mathbb{P}^{n}$, which has
$c_{j}=1$, the standard generator of $H^{2j}(\C\mathbb{P}^{n})$. 
\end{proof}

\subsection{Pontryagin classes}

Since the $j^{th}$ Pontryagin class $P_{j}$ of a real, rank-$k$
vector bundle $\xi$ is just the $2j^{th}$ Chern class of $\xi\otimes\C$
\cite{Milnor-Stasheff}, the form $\Phi P_{j}(\omega)$ will be defined
on the Stiefel bundle $B=V_{k-2j+1}(\xi\otimes\C)$ of complex $(k-2j+1)$-frames
of $\xi\otimes\C$ as the form $\Phi c_{2j}(\omega_{\C})$, using
the natural extension of the connection $\omega$ to $F(\xi\otimes\C).$
However, there is an interpretation of the Pontryagin classes in at
least one case, $P_{1}$, which is independent of a complexification
of $\xi$. 

Let $\xi$ be a real, oriented, rank-$4$ vector bundle over $M$.
Then, $\xi$ is a rank-one, quaternionic vector bundle over $M$,
that is, the bundle $F(\xi)$ of oriented frames of $\xi$ defines
a bundle of (non-integrable) quaternionic structures on $\xi$. Each
frame, that is, each quaternionic structure of $\xi_{x}$, determines
3 complex structure tensors on $\xi_{x}$, $I$, $J$, and $K$, with
$IK=-J$. $I$ is defined by $I(e_{1})=e_{2}$, $I(e_{3})=e_{4}$,
and of course $I(e_{2})=-e_{1}$and $I(e_{4})=-e_{3}$. Similarly,
$K$ is defined at the same frame by $K(e_{1})=e_{3}$ and $K(e_{2})=e_{4}$,
and $J$ is defined by $J=KI$. This bundle of quaternionic structures
produces 2 complementary bundles of complex structures (corresponding
to $I$ and $K$, to be specific). If $H_{1}$ and $H_{2}$ are the
two subgroups of $SO(4)$ corresponding to the complex-linear automorphisms
with complex determinant 1 with respect to $I$ and $K$, respectively,
each being a representation of $SU(2)$ in $SO(4)$, then the associated
bundles $B_{1}:=F(\xi)\times_{SO(4)}SO(4)/H_{1}$ and $B_{2}:=F(\xi)\times_{SO(4)}SO(4)/H_{2}$
are each $\R\mathbb{P}^{3}$-bundles over $M$.  In general, of course,
the complex-structure tensors $I$ and $K$ will not be well-defined
on all of $M$, but if so, such as for the tangent bundle of a hyperk\"{a}hler
manifold of real dimension 4, then they would give two dual complex
structures on $\xi$. Such complex structures would correspond to
global sections of $B_{1}$ and $B_{2}$.

Now, it will not be the case that, for either bundle, the term $P_{1}(\Psi,\Psi)=0$
as was the case in the previous situations. However, since \begin{eqnarray*}
P_{1}(\Omega,\Omega) & = & \lambda Tr(\Omega\wedge\Omega)\\
 & = & -\lambda(\Omega_{12}\Omega_{12}+\cdots+\Omega_{34}\Omega_{34})\\
 & = & -\lambda\sum_{i<j}\Omega_{ij}\wedge\Omega_{ij}\end{eqnarray*}
 for $\Omega\in so(4)$, and since, in this situation, the decomposition
$so(4)=\frak{h}_{1}\oplus\frak{p}_{1}=\frak{p}_{2}\oplus\frak{h}_{2}=su(2)\oplus su(2)$
splits as Lie algebras rather than just as a reductive complement,
then, for either $B_{1}$ or $B_{2}$ the decomposition $\omega=\phi_{i}+\psi_{i}$
of the connection satisfies $[\phi_{i},\phi_{i}]\subset\frak{p}_{i}$
(that is, the form takes values in $\frak{p}_{i}$) and $[\psi_{i},\psi_{i}]\subset\frak{h}_{i}$.
Thus, by the decomposition of $\Psi_{i}=\Omega_{\frak{h}_{i}}+\frac{1}{2}[\phi_{i},\phi_{i}]_{\frak{h}_{i}}=\Omega_{\frak{h}_{i}}$,
so $P_{1}(\Psi_{i},\Psi_{i})=P_{1}(\Omega_{\frak{h}_{i}},\Omega_{\frak{h}_{i}})$,
denoting by $\Psi_{i}$, $i=1,2$, the corresponding curvature forms,
for $B_{1}$and $B_{2}$ (both of which can be viewed as forms in
$F(\xi)$ having values in $so(4)$), then \[
P_{1}(\Psi_{1},\Psi_{1})+P_{1}(\Psi_{2},\Psi_{2})=P_{1}(\Omega,\Omega),\]
so that while neither one of the transgressive forms has differential
the Pontryagin class, their sum does, \[
d\Phi P_{1}(\omega_{1})+d\Phi P_{1}(\omega_{2})=P_{1}(\Omega),\]
where of course $\omega_{1}$ and $\omega_{2}$ refer to the two distinct
decompositions of the bundle of frames into associated bundles (even
though $\omega$ is the same in both cases). 

\begin{prop}
Let $\xi$ be a real, oriented, rank-$4$ vector bundle over a compact
manifold $M$. Let $B_{1}$ and $B_{2}$ be given by $B_{1}:=F(\xi)\times_{SO(4)}SO(4)/H_{1}$
and $B_{2}:=F(\xi)\times_{SO(4)}SO(4)/H_{2}$ as above. Let $\alpha$
be a smooth singular $4$-chain in $M$. Choose generic sections $\sigma_{1}$
of $B_{1}$ and $\sigma_{2}$ of $B_{2}$ with a finite set of singular
points $\{ p_{j}\}$, singular for either $\sigma_{1}$ or $\sigma_{2}$,
or both, interior to $4$-simplices in $\alpha$, with nondegenerate
singularities of indices $a_{1l}$ and $a_{2l}$ at $p_{l}$. Then,
\begin{eqnarray*}
\int_{\alpha}P_{1}(\Omega) & = & \sum_{l=1}^{m}a_{1l}+a_{2l}+\int_{\partial\alpha}\sigma_{1}^{*}(d\Phi P_{1}(\omega_{1}))+\sigma_{2}^{*}(d\Phi P_{1}(\omega_{2})).\end{eqnarray*}

\end{prop}
\begin{proof}
Let $\alpha$ be a smooth singular $4$-cycle in $M$. $B_{1}$ and
$B_{2}$ give rise to two $\R^{4}$-bundles on $M$, which as above
will have generic sections with a discrete set of nondegenerate zeros
on $\alpha$, corresponding to sections $\sigma_{1}$ of $B_{1}$
and $\sigma_{2}$ of $B_{2}$ with a finite set of singular points
$\{ p_{j}\}$, singular for either $\sigma_{1}$ or $\sigma_{2}$,
or both, interior to $4$-simplices in $\alpha$, which are limits
of maps $\sigma_{i}:\partial B(p_{j},\epsilon)\rightarrow\R\mathbb{P}^{3}$
of finite degree when lifted to $\widetilde{\sigma}_{i}:\partial B(p_{j},\epsilon)\rightarrow S^{3}$,
$i=1,2$, that is, for $\epsilon>0$ sufficiently small, $\left(\left.\widetilde{\sigma}_{i}\right|_{\partial B(p_{j},\epsilon)}\right)_{*}:H_{3}(\partial B(p_{j},\epsilon))\rightarrow H_{3}(\widetilde{\pi^{-1}(p_{j})})$
given by $[\partial B(p_{j},\epsilon)]\mapsto a_{i,j}[\widetilde{\pi^{-1}(p_{j})}]$,
$a_{ij}\in\Z$. As above, with the projections $\pi_{1}:B_{1}\rightarrow M$
and $\pi_{2}:B_{2}\rightarrow M$, \begin{eqnarray*}
\int_{\alpha}P_{1}(\Omega) & = & \lim_{\epsilon\downarrow0}\int_{\alpha\backslash B(p_{1},\epsilon)\cup\cdots\cup B(p_{m},\epsilon)}\sigma_{1}^{*}(d\Phi P_{1}(\omega_{1}))+\sigma_{2}^{*}(d\Phi P_{1}(\omega_{2}))\\
 & = & \lim_{\epsilon\downarrow0}\sum_{l=1}^{m}\int_{(\sigma_{1})_{*}(\partial B(p_{l},\epsilon))}\Phi P_{1}(\omega_{1})+\int_{(\sigma_{2})_{*}(\partial B(p_{l},\epsilon))}\Phi P_{1}(\omega_{2})\\
 &  & +\int_{\partial\alpha}\sigma_{1}^{*}(d\Phi P_{1}(\omega_{1}))+\sigma_{2}^{*}(d\Phi P_{1}(\omega_{2}))\\
 & = & \sum_{l=1}^{m}a_{1l}\int_{\pi_{1}^{-1}(p_{l})}\Phi P_{1}(\omega_{1})+a_{2l}\int_{\pi_{2}^{-1}(p_{l})}\Phi P_{1}(\omega_{2})\\
 &  & +\int_{\partial\alpha}\sigma_{1}^{*}(d\Phi P_{1}(\omega_{1}))+\sigma_{2}^{*}(d\Phi P_{1}(\omega_{2}))\\
 & = & \sum_{l=1}^{m}a_{1l}\int_{\pi_{1}^{-1}(p_{l})}\sum_{i=0}^{1}A_{i,1-i}P_{1}(\phi_{1},[\phi_{1},\phi_{1}]^{i},\Psi_{1}^{1-i})+\\
 &  & +a_{2l}\int_{\pi_{2}^{-1}(p_{l})}\sum_{i=0}^{1}A_{i,1-i}P_{1}(\phi_{2},[\phi_{2},\phi_{2}]^{i},\Psi_{2}^{1-i})\\
 &  & +\int_{\partial\alpha}\sigma_{1}^{*}(d\Phi P_{1}(\omega_{1}))+\sigma_{2}^{*}(d\Phi P_{1}(\omega_{2}))\\
 & = & \sum_{l=1}^{m}a_{1l}\int_{\pi_{1}^{-1}(p_{l})}A_{0,1}P_{1}(\phi_{1},[\phi_{1},\phi_{1}])+a_{2l}\int_{\pi_{2}^{-1}(p_{l})}A_{0,1}P_{1}(\phi_{2},[\phi_{2},\phi_{2}])\\
 &  & +\sum_{l=1}^{m}a_{1l}\int_{\pi_{1}^{-1}(p_{l})}A_{1,0}P_{1}(\phi_{1},\Omega_{\frak{h}_{1}})+a_{2l}\int_{\pi_{2}^{-1}(p_{l})}A_{1,0}P_{1}(\phi_{1},\Omega_{\frak{h}_{2}})\\
 &  & +\int_{\partial\alpha}\sigma_{1}^{*}(d\Phi P_{1}(\omega_{1}))+\sigma_{2}^{*}(d\Phi P_{1}(\omega_{2}))\\
 & = & \sum_{l=1}^{m}a_{1l}\int_{\pi_{1}^{-1}(p_{l})}P_{1}(\phi_{1},[\phi_{1},\phi_{1}])+a_{2l}\int_{\pi_{2}^{-1}(p_{l})}P_{1}(\phi_{2},[\phi_{2},\phi_{2}])\\
 &  & +0\\
 &  & +\int_{\partial\alpha}\sigma_{1}^{*}(d\Phi P_{1}(\omega_{1}))+\sigma_{2}^{*}(d\Phi P_{1}(\omega_{2})),\end{eqnarray*}
since the integration is over $\pi_{i}^{-1}(p_{l})$, and $\Omega$
is $\pi_{i}$-horizontal. 

Since the form $P_{1}(\omega,[\omega,\omega])$ has integral periods
and generates the transgressive first Pontryagin form of $H^{3}(SO(4),\R)=H^{3}(SU(2),\R)\oplus H^{3}(SU(2),\R)$,
the projection $SO(4)\rightarrow SO(4)/H_{1}$ pulls the generator
of $H^{3}(SO(4)/H_{1},\R)$, which is $P_{1}(\phi_{1},[\phi_{1},\phi_{1}])$,
back to $P_{1}(\omega,[\omega,\omega])$, and so \[
\int_{\pi_{1}^{-1}(p_{l})}P_{1}(\phi_{1},[\phi_{1},\phi_{1}])=1;\]
similarly with the other projection as well. Thus\begin{eqnarray*}
\int_{\alpha}P_{1}(\Omega) & = & \sum_{l=1}^{m}a_{1l}+a_{2l}+\int_{\partial\alpha}\sigma_{1}^{*}(d\Phi P_{1}(\omega_{1}))+\sigma_{2}^{*}(d\Phi P_{1}(\omega_{2})).\end{eqnarray*}

\end{proof}
\begin{example}
As an example of this decomposition, let $\xi$ be the tangent bundle
to $S^{4}$. Since $S^{4}=\mathbb{HP}^{1}$, $S^{4}$ admits a global
quaternionic structure, though it admits no global almost-complex
structures. However, if $p$ is the South pole, $p=(-1,0,0,0,0)$,
on $S^{4}\backslash\{ p\}$, there are certainly global almost-complex
structures. The standard complex structures can be described by parallel
transport of a given pair of complex structures at the North pole
$n:=-p$, corresponding to the canonical frame, along longitudes.
For these sections $\sigma_{1}$ of $B_{1}$ and $\sigma_{2}$ of
$B_{2}$, there is only one singular point, for both sections, at
$p$. The section $\left.\sigma_{1}\right|_{x_{1}=-\sqrt{1-\epsilon^{2}}}:S_{\epsilon}^{3}\rightarrow\left.B_{1}\right|_{S_{\epsilon}^{3}}\cong S_{\epsilon}^{3}\times\mathbb{RP}^{3}$,
as a map $\sigma_{1}S^{3}\rightarrow\mathbb{RP}^{3}$, lifts to a
map $\widetilde{\sigma_{1}}:S_{\epsilon}^{3}\rightarrow S^{3}$ of
degree $2$, and $\sigma_{2}$ similarly lifts to a map of degree
$-2$. 
\end{example}
\begin{rem}
An eventual goal of these relative classes would be to construct a
combinatorial procedure to determine the Pontryagin classes of a closed
manifold, since it is well-known that they are topological invariants
\cite{Novikov}. It is possible to begin the procedure for an oriented
$4$-manifold $M$ based on this result, but the details are not apparent.
Given such a manifold $M^{4}$, each $4$-simplex of a fixed triangulation
of $M$ admits a standard hyperk\"{a}hler structure, with prescribed
behavior at the boundary (given as that of the boundary of the 4-ball
in quaternionic 1-space), so that the computation reduces to the 3-skeleton.
As with the Euler class, the boundary terms cancel on the interiors
of the 3-simplices due to reversal of orientation of the sections,
so the remaining calculations should lie on the 2-skeleton. 
\end{rem}

\section{Secondary characteristic classes}

One of the most extensive uses which has been made with the construction
of \cite{CS} has been the construction of secondary characteristic
classes. If $P$ is a polynomial of degree $k$ so that $P(\Omega)$
is integral for all $\Omega$, that is, if it has integral periods,
then when $\omega$ is a connection for which $P(\Omega)=0$ (as a
form, not just as a cohomology class), the Chern-Simons transgression
$TP(\omega)$ will be closed, generating a cohomology class in $H^{2k-1}(E,\R)$.
Of more interest is the construction, from that class, of a cohomology
class in the base $M$ modulo integral classes. In the case of a principal
bundle the existence of such a class follows by passing to the universal
bundle, where every cocycle is a coboundary on the total space, and
so the mod-$\Z$ reduction of $TP(\omega)$ will be a lift of a cocycle
on the base. 

The forms $\Phi P(\omega)$ can, in some cases, be more directly seen
to be lifts, using the obstruction information determined by the characteristic
class $P(\Omega)$. Note that the method of proof used by Chern and
Simons will not work in this situation, and that the forms $\Phi P(\omega)$,
and so the secondary characteristic classes determined by them, are
not the same as the Chern-Simons classes. 

\begin{thm}
If $\xi$ is either a rank-$2k$ real oriented vector bundle, or a
rank-$k$ complex vector bundle, over a compact manifold $M$, and
if, respectively, the form $e(\Omega)=0$ (resp., $c_{j}(\Omega)=0$
for some $j$), then the corresponding form $\Phi e(\omega)$ (resp.,
$\Phi c_{j}(\omega)$) is well-defined as an element of $H^{2k-1}(M,\R/\Z)$
(resp., $H^{2j-1}(M,\R/\Z)$.
\end{thm}
\begin{proof}
For the Euler class of an oriented rank-$2k$ vector bundle, where
$B$ is the sphere bundle: if $e(\Omega)=0$ as a form, of course
$\Phi e(\omega)$ will be closed, and for any section $\sigma$ over
the $2k$-skeleton, $\sigma^{*}(\Phi e(\omega))$ will be a closed
form on $M$ which lifts to $\Phi e(\omega)$, defining a secondary
characteristic class on $M$ modulo the choice of section $\sigma$. 

The Gysin sequence of the $(2k-1)$-sphere bundle $\pi:B\rightarrow M$,
\[
\cdots\begin{array}{c}
\pi_{*}\\
\rightarrow\\
\,\end{array}H^{r}(M,R)\begin{array}{c}
e\wedge\\
\rightarrow\\
\,\end{array}H^{r+2k}(M,R)\begin{array}{c}
\pi^{*}\\
\rightarrow\\
\,\end{array}H^{r+2k}(B,R)\begin{array}{c}
\pi_{*}\\
\rightarrow\\
\,\end{array}H^{r+1}(M,R)\begin{array}{c}
e\wedge\\
\rightarrow\\
\,\end{array},\]
for $r=-1$ yields the split short exact sequence (for any section
$\sigma$)\[
0\begin{array}{c}
\\\rightarrow\\
\,\end{array}H^{2k-1}(M,R)\begin{array}{c}
\pi^{*}\\
\rightarrow\\
\leftarrow\\
\sigma^{*}\end{array}H^{2k-1}(B,R)\begin{array}{c}
\pi_{*}\\
\rightarrow\\
\,\end{array}R.\]
The map $\pi_{*}$ is integration over the fiber \cite[p. 178]{Bott-tu},
so in the case that $R=\R$ the image $\sigma^{*}(\Phi e(\omega))$
is well-defined modulo $\pi_{*}(\Phi e(\omega))=[f]$, where $f(x)=\int_{\pi^{-1}(x)}\Phi e(\omega)$
is an integer-valued (hence constant, since it is continuous) function
by the fact that the characteristic class $e(\Omega)$ has integral
periods. Then, with the coefficient ring $R=\R/\Z$, $\Phi e(\omega)\in ker(\pi_{*})$,
thus there is a unique $U\in H^{2k-1}(M,\R/\Z)$ so that $\pi^{*}(U)=\Phi e(\omega)$,
and since $\pi\sigma=1$, $U=\sigma^{*}\pi^{*}(U)=\sigma^{*}(\Phi e(\omega))$,
independent of choice of $\sigma$. 

Working directly with real cohomology classes, the images $\sigma^{*}(\Phi e(\omega))$
will depend upon the choice of the section $\sigma$, that is, on
a specific splitting of the short exact sequence \[
0\begin{array}{c}
\\\rightarrow\\
\,\end{array}H^{2k-1}(M,\R)\begin{array}{c}
\pi^{*}\\
\rightarrow\\
\leftarrow\\
\sigma^{*}\end{array}H^{2k-1}(B,\R)\begin{array}{c}
\pi_{*}\\
\rightarrow\\
\,\end{array}\R\rightarrow0.\]
These splittings are determined by $Hom(\R,H^{2k-1}(M,\R))$, and
those coming from sections also satisfy an integrality condition since
they also split the integral cohomology sequence\[
0\begin{array}{c}
\\\rightarrow\\
\,\end{array}H^{2k-1}(M,\Z)\begin{array}{c}
\pi^{*}\\
\rightarrow\\
\leftarrow\\
\sigma^{*}\end{array}H^{2k-1}(B,\Z)\begin{array}{c}
\pi_{*}\\
\rightarrow\\
\,\end{array}\Z\rightarrow0,\]
and so such splittings are determined by $Hom(\Z,H^{2k-1}(M,\Z))$,
the free part of $H^{2k-1}(M,\Z)$, again showing that $\Phi e(\omega)$
is well defined on $M$ modulo forms with integral periods, that is,
as an element of $H^{2k-1}(M,\R/\Z)$. 

In the case of a rank-$k$ complex vector bundle, if the form $c_{j}(\Omega)=0$
then, as for the Euler class, there will be a section $\sigma:M^{(2j)}\rightarrow V_{k-j+1}(\xi)$
of the Stiefel bundle of unitary $(k-j+1)$-frames of $\xi$ over
the $2j$-skeleton of $M$. The Stiefel bundle splits as a tower of
sphere bundles \[
V_{k-j+1}(\xi)\rightarrow\cdots\rightarrow V_{2}(\xi)\cong S(\pi^{*}(v^{\perp}))\rightarrow S(\xi)\rightarrow M,\]
where $V_{2}(\xi)$ is the sphere bundle of the orthogonal complement
bundle $v^{\perp}\rightarrow S(\xi)$ with fiber over $v\in S(\xi)$
the orthogonal complement of $\{ v,iv\}$ in the fiber $\xi_{\pi(v)}$.
The fiber at each stage is $S^{2k-2l-1}$, $l=0,\ldots,k-j$. Applying
the Gysin sequence at each stage, with $r=2(j-k)-1$ in the first
stage through $r=-1$ at the last, gives \[
H^{2j-1}(M,R)\cong H^{2j-1}(S(\xi),R)\cong\cdots\cong H^{2j-1}(V_{k-j}(\xi),R),\]
and a split exact sequence (using $H^{2j-1}(V_{k-j}(\xi),R)\cong H^{2j-1}(M,R)$)
for any section $\sigma$ over the $2j$-skeleton of $M$:\[
0\begin{array}{c}
\\\rightarrow\\
\,\end{array}H^{2j-1}(M,R)\begin{array}{c}
\pi^{*}\\
\rightarrow\\
\leftarrow\\
\sigma^{*}\end{array}H^{2j-1}(V_{k-j+1}(\xi),R)\begin{array}{c}
\pi_{*}\\
\rightarrow\\
\,\end{array}R.\]
The proof then proceeds as in the first case, noting that $\pi_{*}$
can still be viewed as integration over the fiber, but over the fiber
of the $S^{2j-1}$-bundle $V_{k-j+1}(\xi)\rightarrow V_{k-j}(\xi)$.
\end{proof}

\end{document}